\documentclass[11pt]{amsart}
\usepackage{amsfonts}
\usepackage{amssymb}
\usepackage[latin1]{inputenc}
\newtheorem{thm}{Theorem}[section]

\newtheorem{lemma}[thm]{Lemma}

\newtheorem{prop}[thm]{Proposition}
\theoremstyle{definition}
\newtheorem{defn}[thm]{Definition}
\theoremstyle{remark}
\newtheorem{rem}[thm]{Remark}
\newtheorem{remark}[thm]{Remark}
\newtheorem{example}[thm]{Example}

\newcommand{\bbR}{\mathbb R}
\newcommand{\bbF}{\mathbb F}

\newcommand{\fg}{\mathfrak g}

\numberwithin{equation}{section}

\newcommand{\rank}{{\rm rank\ }}

\begin{document}

\title[Equivariant normal forms of Poisson structures]{A note on equivariant normal forms of Poisson structures}
\author{Eva Miranda}
\address{Laboratoire Emile Picard, UMR 5580 CNRS, Universit{\'e} Toulouse III}

\email{ miranda@picard.ups-tlse.fr}
\thanks{The first author is supported by Marie Curie EIF postdoctoral fellowship contract number
 EIF2005-024513 and partially supported by the DGICYT project number
BFM2003-03458.}


\author{Nguyen Tien Zung}
\address{Laboratoire Emile Picard, UMR 5580 CNRS, Universit\'{e} Toulouse III}
\email{tienzung@picard.ups-tlse.fr}

\date{Second version, March 16, 2006}
\subjclass{ 53D17}

\keywords{Poisson manifolds;   Linearization;
Normal Forms}

\begin{abstract}
   We prove an equivariant version of the local splitting
  theorem for tame Poisson structures and Poisson actions of compact Lie groups. As a consequence,  we obtain an 
equivariant  linearization result for Poisson structures whose transverse structure has semisimple linear part of 
compact type.

\end{abstract}
\maketitle

\section{Introduction}

The main purpose of this note is to prove an equivariant version of Weinstein's splitting theorem for Poisson
structures \cite{weinstein}. This theorem asserts that in the neighborhood of any point $p$ in a Poisson manifold
$(P^n,\Pi)$ there is a local coordinate system $(x_1,y_1,\dots, x_{2k},y_{2k}, z_1,\dots, z_{n-2k})$ in which the Poisson structure $\Pi$ can be written as
\begin{equation}
\Pi = \sum_{i=1}^k
\frac{\partial}{\partial x_i}\wedge \frac{\partial}{\partial y_i} + \sum_{ij} f_{ij}(z)\frac{\partial}{\partial z_i}\wedge \frac{\partial}{\partial z_j} ,
\end{equation}
where $2k$ is the rank of $\Pi$ at $p$, and $f_{ij}$ are functions
which depend only on the variables $(z_1,\hdots,z_{n-2k})$ and
which vanish at the origin. Geometrically speaking, locally the
Poisson manifold $(P^n,\Pi)$ can be splitted into the direct
product of a $2k$-dimensional symplectic manifold (with the
standard nondegenerate Poisson structure $\Pi_1 = \sum_{i=1}^k
\frac{\partial}{\partial x_i}\wedge \frac{\partial}{\partial
y_i}$) and a $(n-2k)$-dimensional Poisson manifold whose Poisson
structure $\Pi_2 = \sum_{ij} f_{ij}(z)\frac{\partial}{\partial
z_i}\wedge \frac{\partial}{\partial z_j}$ vanishes at the origin.
We want to show that if there is a (local) action of a compact Lie
group $G$ on $P^n$ which fixes the point $p$ and which preserves
$\Pi$, then this splitting can be made equivariantly.

In the special case when $\Pi$ is nondegenerate at $p$ (i.e., $2k = n$), one recovers from Weinstein's theorem the classical Darboux theorem about the local existence of canonical (Darboux) coordinates for symplectic manifolds. We know two methods for proving Darboux theorem: 1) the classical coordinate-by-coordinate construction method; and 2) the path method due to Moser \cite{moser}. Weinstein's proof of the splitting theorem \cite{weinstein} is also based on the first method (coordinate by coordinate construction). However, this classical method does not seem to work in the equivariant situation, while the path method can be used to prove the equivariant Darboux theorem \cite{wei1}.

In the same spirit, we will try to use the path method to prove an
equivariant version of the splitting theorem for Poisson structures.
In doing so, we encounter a technical condition, which we call the
\emph{tameness condition}: a smooth Poisson structure $\Pi$ on a
manifold $P^n$ is called \emph{tame} if for any two smooth Poisson
vector fields $X, Y$ on $P^n$ (which may depend on some parameters)
which are tangent to the symplectic leaves the function
$\Pi^{-1}(X,Y)$ is smooth (and depends smoothly on the parameters).
We will devote Section 2 of this note to the tameness condition, in
order to convince the reader that it is an interesting condition,
and many ``reasonable'' Poisson structures satisfy it. For example,
if the linear part of the transverse Poisson structure at a point
$p$ has semisimple  type, then the Poisson structure is tame near
$p$.

Now we can formulate the main result of this note:

\begin{thm} \label{thm:main}
Let $(P^n,\Pi)$ be a smooth Poisson manifold, $p$ a point of $P$, $2k = \rank \Pi (p)$, and $G$ a compact Lie group which acts on $P$ in such a way that the action preserves $\Pi$ and fixes the point $p$. Assume that the Poisson structure $\Pi$ is tame at $p$. Then there is a smooth canonical local coordinate system $(x_1,y_1,\dots, x_{2k},y_{2k}, z_1,\dots, z_{n-2k})$ near $p$, in which the Poisson structure $\Pi$ can be written as
\begin{equation}
\Pi = \sum_{i=1}^k
\frac{\partial}{\partial x_i}\wedge \frac{\partial}{\partial y_i} + \sum_{ij} f_{ij}(z)\frac{\partial}{\partial z_i}\wedge \frac{\partial}{\partial z_j},
\end{equation}
and in which the action of $G$ is linear and preserves the subspaces $\{x_1 = y_1 = \hdots x_k = y_k = 0\}$ and $\{z_1 = \hdots = z_{n-2k} = 0\}$.
\end{thm}

\begin{remark}
i) We do not know whether the tameness condition is really necessary, or if it is because our method is not good enough. We notice that this condition is also implicitly present in the papers of Ginzburg and Weinstein \cite{ginzburgweinstein} and of Aleekseev and Meinrenken \cite{alek}, \cite{alek2}, which involve the path method in Poisson geometry. \\
ii) The above theorem also holds in the analytic (i.e., real analytic or holomorphic) setting, with basically the same proof. The analytic version of this equivariant theorem is used by Philippe Monnier and the second author  in their study of normal forms of vector fields on Poisson manifolds \cite{zungmonnier}. We hope that our result can be useful in the study of equivariant Hamiltonian systems as well. \\
iii) If the action of $G$ on $(P^n,\Pi)$ is Hamiltonian (with an equivariant momentum map), then there is another approach to this equivariant splitting problem, based on the Nash-Moser method, which does not need the tameness condition. We will consider this issue in a separate work.
\end{remark}

The above theorem will be proved in Section 3 of this note. In
Section 4 we will combine this theorem with linearization results
of Conn \cite{conn} and Ginzburg \cite{ginzburg} to obtain an
equivariant linearization theorem (see Theorem
\ref{thm:EquivLin}).

\vspace{5mm}

{\bf Acknowledgements.} We are indebted to Mich{\`e}le Vergne for
drawing our attention to the paper of Dixmier \cite{dixmier} and
for pointing out its relation to the  division property stated in
section 2.

We would like to thank Viktor Ginzburg for his useful comments and
suggestions on the problem. We would also like to thank David
Mart{\'\i}nez-Torres for carefully reading a previous version of this
preprint and pointing out some misprints.

\section{Tame Poisson structures}

We will denote by $\Pi^{-1}$ the covariant tensor dual to the
Poisson tensor $\Pi$ of a Poisson manifold $(P^n,\Pi)$, i.e. the
symplectic form on symplectic leaves.  If $X,Y$ are vector fields
on $P^n$ which are tangent to the symplectic leaves, then
$\Pi^{-1}(X,Y)$ is well-defined. In particular, if $X = X_h$ is
the Hamiltonian vector field of a function $h$ on $(P^n,\Pi)$ then
$\Pi^{-1}(X,Y) = - Y(h)$. Recall that a Poisson vector field is a
vector field which preserves the Poisson structure.

\begin{defn}Let $(P^n,\Pi)$ be a smooth Poisson manifold and $p$ a point in $P$. We will say that $\Pi$ is \emph{tame}
at $p$ if for any pair $X_t,Y_t$  of germs of smooth Poisson
vector fields near $p$ which are tangent to the symplectic
foliation of $(P^n,\Pi)$ and which may depend smoothly on a
(multi-dimensional) parameter $t$, then then the function
$\Pi^{-1}(X_t,Y_t)$ is smooth and depends smoothly on $t$.
\end{defn}

The tameness condition is a kind of homological condition. In particular, if the parametrized germified first Poisson cohomology group, which we will denote by  $H^1_\Pi(P^n,p)$, vanishes, then $\Pi$ is tame at $p$. Indeed, $H^1_\Pi(P^n,p) = 0$ means that if $X_t$ is a germ of Poisson vector field near $p$ which depends smoothly on a parameter $t$, then we can write
$X_t = X_{h_t}$ where $h_t$ is a germ of smooth function near $p$ which depends smoothly on the parameter $t$. Hence
$\Pi^{-1}(X_t,Y_t) = - Y_t(h_t)$ is smooth.

In particular, it is known that if $\fg$ is a compact semi-simple Lie algebra, and $(\fg^*,\Pi_{lin})$ is the dual of $\fg$ equipped with the corresponding linear Poisson structure then $H^1_{\Pi_{lin}}(\fg^*,0) = 0$ (see \cite{conn}). Hence our first example of tame Poisson structures:

\begin{example}
Any smooth Poisson structure $\Pi$, which vanishes at a point $p$
and whose linear part at $p$ corresponds to a compact semisimple
Lie algebra $\fg$, is tame at $p$. Indeed, in this case, according
to Conn's smooth linearization theorem \cite{conn}, $(P^n, \Pi)$
is locally isomorphic near $p$ to $(\fg^*,\Pi_{lin})$, and
therefore $H^1_\Pi(P^n,p) = 0$.
\end{example}

If $X$ is not Hamiltonian (and maybe not even Poisson) but can be written as $X = \sum_{i=1}^m f_i X_{g_i}$ where $f_i, g_i$ are smooth functions, then $\Pi^{-1}(X,Y) = - \sum_{i=1}^m f_i Y(g_i)$ is still smooth. This leads us to:

\begin{defn}
We say that a smooth (resp real analytic) Poisson structure $\Pi$
satisfies the \emph{ smooth division property} (resp
\emph{analytic division property}) at a point $p$ if the
Hamiltonian vector fields generate the space of vector fields
tangent to the associated symplectic foliation near $p$. More
precisely, for any germ of smooth (resp. analytic) vector field
$Z$ -which may depend smoothly (resp. analytically) on some
parameters- which is tangent to the symplectic foliation there
exists a finite number of germs of smooth (resp. analytic)
functions $f_1,\hdots,f_m,g_1,\hdots,g_m$ -which depend smoothly
(resp. analytically) on the same parameters as $Z$- such that $Z =
\sum f_i X_{g_i}$.
\end{defn}

Clearly, if $\Pi$ satisfies the division property at a point $p$,
then it is tame at $p$. A natural question is to know  which Poisson
structures satisfy the division property. In particular, is it true
that all linear Poisson structures satisfy the division property at
the origin? In the appendix we prove that low-dimensional Lie
algebras satisfy the division property at the origin. Namely

\begin{prop} \label{prop:DivisionDim3}
Any linear Poisson structure in dimension 2 or 3 has the division
property at the origin.
\end{prop}

In the higher-dimensional case, a result of Dixmier \cite{dixmier}
says (in our language) that if $\Pi$ is a linear Poisson structure
which corresponds to a semisimple Lie algebra then it has the
analytic division property at the origin (mainly due to the fact
that the singular set has codimension $3$ in this case). We would
conjecture that Dixmier's result also holds in the smooth case. On
the other hand, one can probably produce linear Poisson (non
semisimple) structures which do not satisfy the division property
(similiar to Dixmier's counterexample $3.3$ in \cite{dixmier}).

It is not difficult to construct examples of Poisson structures with
a trivial 1-jet which are not tame.

\begin{example}
Consider the Poisson structure
$\Pi=x^4 \frac{\partial}{\partial x}\wedge \frac{\partial}{\partial y}$ on $\bbR^2$.
The  following vector fields are Poisson and  tangent to the symplectic foliation:
 $$X=x^2  \frac{\partial}{\partial x}+2xy  \frac{\partial}{\partial y}, \quad Y=x  
\frac{\partial}{\partial y},$$
but $\displaystyle \Pi^{-1}(X,Y)=\frac{1}{x}$ is not smooth at the origin. So this Poisson structure is not tame.
\end{example}

Recall that if  $\Pi = \sum_{i=1}^k
\frac{\partial}{\partial x_i}\wedge \frac{\partial}{\partial y_i} + \sum_{ij} f_{ij}(z)\frac{\partial}{\partial z_i}\wedge \frac{\partial}{\partial z_j}$ in a local canonical coordinate system in the neighborhood of a point $p$, then
$\Pi_2 = \sum_{ij} f_{ij}(z)\frac{\partial}{\partial z_i}\wedge \frac{\partial}{\partial z_j}$ is called the \emph{transverse} Poisson structure of $\Pi$ at $\Pi$. Up to local Poisson isomorphisms, this Poisson transevrse structure is unique, i.e. it does not depend on the choice of local canonical coordinates, see, e.g., \cite{dufourzung,weinstein}. The following lemma shows that, to verify the tameness condition, it is sufficient to check it in the transverse direction to the symplectic leaf:

\begin{lemma}\label{lemma:transverse}  A smooth Poisson structure $\Pi$
is tame at a point $p$ if and only if the transverse Poisson structure of $\Pi$ at $p$ is tame at $p$.
\end{lemma}

\begin{proof}
Write $\Pi = \Pi_1 + \Pi_2 = \sum_{i=1}^k
\frac{\partial}{\partial x_i}\wedge \frac{\partial}{\partial y_i} + \sum_{ij} f_{ij}(z)\frac{\partial}{\partial z_i}\wedge \frac{\partial}{\partial z_j}$ in a local canonical coordinate system near $p$. For each germ of vector field $X$ near $p$ write $X = X_{hor} + X_{vert}$, where $X_{hor}$ is the ``horizontal part'' of $X$, i.e. is a combination of the vector fields
$\frac{\partial}{\partial x_i}, \frac{\partial}{\partial y_i}$, and $X_{vert}$ is the ``vertical part'' of $X$, i.e. is a combination of the vector fields $\frac{\partial}{\partial z_i}$. If $X$ is a smooth Poisson vector field for $\Pi$, then $X_{hor}$ (resp. $X_{vert}$) may be viewed as a Poisson vector field for $\Pi_1$ (resp., $\Pi_2$) which depends smoothly on parameters $z_i$ (resp., $x_i,y_i$). We have $ \Pi^{-1}(X,Y)= \Pi_1^{-1}(X_{hor},Y_{hor})+\Pi_2^{-1}(X_{vert},Y_{vert})$.
The term $\Pi_1^{-1}(X_{hor},Y_{hor})$ is always smooth (provided that $X$ and $Y$ are smooth), and so the smoothness of
$ \Pi^{-1}(X,Y)$ is equivalent to the smoothness of $\Pi_2^{-1}(X_{vert},Y_{vert})$. The lemma then follows easily.
\end{proof}

\section{Proof of the equivariant splitting theorem}

In this section we will give a proof of Theorem  \ref{thm:main}.
It uses coupling tensors for Poisson manifolds, so we will first
recall a result of Yu. Vorobiev  about coupling tensors (see,
e.g., \cite{dufourzung,voro}). The proof of the theorem consists
of three steps. In the first step we prove that we can assume that
the action of our compact Lie group $G$ is linear and that the
symplectic foliation is normalized (i.e. is the same as in the
splitting theorem). In the second step we construct a path of
$G$-invariant  Poisson structures connecting the initial Poisson
structure to the splitted one. Finally, in the last step, we use
this path of Poisson structures and the averaging method to
construct a flow which intertwines with the action of $G$ and
whose time-1 map moves the initial Poisson structure to the
splitted one, thus giving an equivariant splitting of our Poisson
structure.

\subsection{Preliminaries: coupling tensors} \hfill

Let $\pi: E \longrightarrow S$ be a submersion over a manifold $S$
and let $T_V E=\ker d\pi$. An Ehresmann connection on $E$ is a
splitting of the tangent bundle of $E$ as $TE=T_V E\oplus T_H E$.
 We call
$T_H E$ the horizontal space. Denote by $\mathcal{V} ^1_V(E)$ the
set of vertical vector fields.
 We can associate to this splitting a $\mathcal{V}^1_V(E)$-valued
 $1$-form $\Gamma\in\Omega^1(E)\otimes \mathcal{V}^1_V(E)$ such that
 $\Gamma(Z)=Z$ for any vertical vector field. Then the horizontal space can be 
written as
 $T_H E=\{X\in TE,\quad \Gamma(X)=0 \}$.
  We can define
 the  horizontal lifting of vector fields from $S$ to $E$.
 In the same way, we may associate a parallel transport to $\Gamma$
 which is smooth, a curvature form and a covariant derivative
 (for details see for example  \cite{dufourzung}).

Consider now the case when $S$ is a symplectic leaf of a  Poisson manifold
$(P,\Pi)$. We can consider a neighbourhood $E$ of $S$ and  submersion $\pi: E \longrightarrow S$
whose restriction to $S$ is the identity.

There is a natural smooth Ehresmann connection where the
horizontal subbundle is spanned by the Hamiltonian vector fields
$X_{f\circ\pi}$.

 We can also
associate to it a $2$-form $\mathbb F\in \Omega^2(S)\otimes
{\mathcal C}^{\infty}(E)$ defined as $$\mathbb
F(X_{f\circ\pi},X_{g\circ\pi})=\langle
\Pi,\pi^*{df}\wedge\pi^*{dg}\rangle.$$

Recall  that  we have an induced  transverse Poisson structure
$\Pi_{Vert}$ on the vertical space.

 The triple
$(\Pi_{Vert},\Gamma,\mathbb F)$ is called the geometric data
associated to the Poisson manifold $(P,\Pi)$ in a neighbourhood of
a symplectic leaf. In \cite{voro}, Vorobjev studies the
reconstruction problem from given geometric data. That is given a
triple of smooth geometric data he gives compatibility conditions
 that guarantee the existence of a Poisson structure with the
 given geometric data. Those compatibility conditions come from the Schouten condition $[\Pi,\Pi]=0$
 imposed on the bivector field $\Pi$ reconstructed from the
 geometric data.

Assume that we are given  $(\Pi_{Vert},\Gamma,\mathbb F)$ on a
fibration $\pi: E \longrightarrow S$, where $\Gamma$  is an
Ehresmann connection on $E$, $\Pi_{Vert}$  a vertical bivector
field, and $\mathbb F \in \Omega^2(S) \otimes \mathcal{C}^\infty(E)$
a nondegenerate $\mathcal{C}^\infty(E)$-valued 2-form on $S$.

We will need the following characterization of
geometrical data which come from a Poisson structure:

\begin{thm}[Vorobiev \cite{voro}]
\label{voro}
The triple  $(\Pi_{Vert},\Gamma,\mathbb F)$ on a
fibration $\pi: E \longrightarrow S$ determines a Poisson structure on
$E$  if and only if $\bbF$ is nondegenerate and the following four
compatibility conditions are satisfied:
\begin{eqnarray}
 &&[\Pi_{Vert} ,\Pi_{Vert} ]=0, \\
&&L_{Hor(u)}(\Pi_{Vert}) =0 \ \ \forall\ u \in \ \nu ^1_V(E), \\
 &&\partial_\Gamma \mathbb F=0, \\
 && Curv_\Gamma(u,v)=\nu^{\sharp}(d(\mathbb F(u,v))) ,
\end{eqnarray}
where $\partial_\Gamma$ stands for the covariant derivative and
$\nu^\sharp$ stands for the map from $T^*E$ to $TE$ defined by
$\langle \nu^\sharp(\alpha), \beta  \rangle = \langle
\nu,\alpha\wedge \beta\rangle. $
\end{thm}

\begin{rem}
We may think of $\Pi$ as the coupling of $\Pi_{Vert}$ with $\mathbb F $ by  $\Gamma$.
  This so-called coupling
 method is a generalization of the minimal coupling procedure
 established for symplectic fibrations by Guillemin, Lerman and
 Sternberg \cite{guillemin}, \cite{sternberg}.
\end{rem}

\subsection{First step of the proof:  linearization of the group action} \hfill

Consider an action $\rho: G \times P^n \rightarrow P^n$ of a compact Lie group $G$ on a
 Poisson manifold $(P^n,\Pi)$, which fixes a point $p \in P^n$ and preserves the Poisson
  structure $\Pi$. Denote by $S$ the local symplectic leaf through $p$. Note that $S$ is invariant under the action of $G$. According to Bochner's theorem \cite{bo}, the action of $G$ is linearizable near $p$, i.e., there is a local coordinate system in which the action is linear. Moreover, we may assume that $S$ is linear in these cordinates. Since linear representations of compact Lie groups are completely reducible, there is a local submanifold $N$ (which is also linear in these coordinates), which is invariant under the action of $G$ and which is transverse to $S$ at $p$. The following lemma says that we can choose this coordinate system in such a way that the symplectic foliation of $(P^n,\Pi)$ will also be the same as in the splitting theorem.

\begin{lemma}\label{samefoliation} With the above notations, there is a local system of coordinates near $p$ in which the action of $G$ is linear, the submanifolds $S$, $N$ are linear, and the local symplectic leaves near $p$ are direct products of $S$ with symplectic leaves of the transverse Poisson structure on $N$.
\end{lemma}

\begin{proof} We can start with a first coordinate system in which the action of $G$ is linear and the submanifolds $S$, $N$ are linear. Denote by $p_1$ the linear projection from a sufficiently small neighborhood $U$ of $p$ in $P^n$ to $S$ which projects $N$ to $p$. Define another (a-priori nonlinear) projection $p_2$, from $U$ to $N$, as follows: Denote by $\Gamma$ the Ehresmann
connection associated to the Poisson structure $\Pi$ and the projection $p_1$.
 For each $x\in U$, let $\alpha_x(t)$ be the linear path joining $p_1(x)$ to the origin $p$ in $S$, with $\alpha_x(0) = p_1(x)$ and $\alpha_x(1) = p$. Denote by $\hat{\alpha}_x$ the horizontal lift of $\alpha_x$
through $x$ with respect to $\Gamma$. Then we take $p_2(x)=\hat{\alpha}_x(1) \in N$.

By construction both projections are smooth and $G$-equivariant: The projection $p_1$ is equivariant since $N$ is $G$-invariant
and $p_2$ is equivariant because the action of $G$ preserves $\Pi$ and therefore the parallel transport is equivariant.

Now consider the $G$-equivariant local diffeomorphism
\begin{equation*}
\begin{array}{ccc}\phi: & U\longrightarrow & S\times N
\\ & x\longmapsto & (p_1(x),p_2(x))
 \end{array}
\end{equation*}

Since the parallel transport preserves the Poisson structure,
$\phi$ takes the Poisson structure on $U$ to a Poisson structure
on $S\times N$ which has as symplectic leaves the product of the
symplectic leaves on $N$ with $S$. This ends the proof of the
lemma.
\end{proof}

\subsection{Second step: constructing a path of Poisson structures} \hfill

After the first step, we can now assume that $P = N \times S$, and
the Poisson structure $\Pi$ has the same symplectic leaves as the
splitted Poisson structure $\tilde{\Pi} =  \Pi_S + \Pi_N$, where
$\Pi_S$ is the standard nondegenerate Poisson structure on $S$ and
$\Pi_N$ is the transverse Poisson structure on $N$, and both $\Pi$
and $\tilde\Pi$ are invariant under our linear action of $G$. We
will assume that $\Pi$ is tame at $p$, or equivalently, the
transverse Poisson structure $\Pi_N$ is tame at the origin.

\begin{lemma}
\label{path}
With the above notations and assumptions, there is a smooth path of $G$-invariant Poisson structures $\Pi_t$, $t \in [0,1]$, on (a neighborhood of the origin in) $N \times S$, such that $\Pi_0 = \Pi$, $\Pi_1 = \Pi_S + \Pi_N$, and which have the same symplectic foliation for all $t \in [0,1]$.
\end{lemma}

\begin{proof}

We denote by $\omega_0$ the symplectic structure induced on the
symplectic leaves by $\Pi_0 = \Pi$. In the same way we denote by
$\omega_1$ the symplectic structure induced  by $\Pi_1 = \Pi_S + \Pi_N$ on the
same symplectic foliation. Consider the linear path of $2$-forms
\begin{equation}
\omega_t=t\omega_1+(1-t)\omega_0 .
\end{equation}
This is a path of smooth closed $2$-forms on each symplectic leaf of the common symplectic foliation. We want to show that, for each $t$ there is a smooth bivector field $\Pi_t$ which corresponds to $\omega_t$. Then, automatically, $\Pi_t$ is a Poisson structure because of the closedness of $\omega_t$, has the same symplectic foliation as $\Pi_0$ and $\Pi_1$, and is $G$-invariant.

Denote by  $(\Pi_N,\Gamma_0, \mathbb F_0)$ and  $(\Pi_N,\Gamma_1,
\mathbb F_1)$, the geometric data associated to the Poisson
structures $\Pi_0 = \Pi$ and $\Pi_1 = \Pi_N + \Pi_S$ with respect to the projection $p_1: N\times S \rightarrow S$ (remark that, by construction, they have the same vertical component, which is equal to $\Pi_N$). We will use Vorobjev's Theorem \ref{voro} to
construct $\Pi_t$ and to prove its smoothness. In other words, we will construct
geometric data $(\Pi_N,\Gamma_t, \mathbb F_t)$, which will be shown to be smooth
and satisfy the compatibility conditions of Theorem \ref{voro}, so they will give
rise to a smooth Poisson structure $\Pi_t$.

\vspace{3mm} {\it Construction and smoothness of
$\Gamma_t$:} \vspace{3mm}

In order to construct the connection $\Gamma_t$, it is enough to show how to lift each vector field $X$ on $S$ horizontally with respect to $\Gamma_t$. The horizontal lift $X_t$ of $X$ with respect to $\Gamma_t$ is uniquely characterized by $\omega_t$ (the would-be associated symplectic form on the symplectic leaves) and by the following two conditions:
\begin{enumerate}
\item The vector field $X_t$ is tangent to the common symplectic foliation of $\Pi_0$ and $\Pi_1$, and its projection to $S$ by $p_1$ is $X$.
\item   $\omega_t(X_t,Z)=0$ for any vertical vector field $Z$.
\end{enumerate}

Denote by $X_0$ and $X_1$ the horizontal lift of $X$ with respect to $\Gamma_0$ and $\Gamma_1$ respectively. We will show that
\begin{equation}
X_t = (1-t) X_0 + tX_1.
\end{equation}
(Then the smothness of $X_t$, and hence of $\Gamma_t$, is
automatic). It is clear that $(1-t) X_0 + tX_1$ is tangent to the
symplectic foliation and projects to $X$ under $p_1$. It remains
to show that
\begin{equation*}
\omega_t ((1-t) X_0 + tX_1, Z) = 0
\end{equation*}
for any vertical vector field $Z$ on $N \times S$. Indeed, denoting $W = X_0 - X_1$, we have
\begin{equation*}
\begin{array}{ccl}
&   & \omega_t((1-t) X_0 + tX_1,Z) \\
& = & t\omega_1((1-t) X_0 + tX_1,Z)+(1-t)\omega_0((1-t) X_0 + tX_1,Z) \\
& = & t\omega_1(X_1+(1-t)W,Z)+(1-t)\omega_0(X_0-tW,Z) \\
& = & t\omega_1(X_1,Z) +(1-t)\omega_0(X_0,Z) + t(1-t) [\omega_1(W,Z) - \omega_0(W,Z)] .
\end{array}
\end{equation*}

Since $X_1$ and $X_0$ are the horizontal lifts of $X$ with respect
to $\Gamma_1$ and $\Gamma_0$, the terms $\omega_1(X_1,Z)$ and
$\omega_0(X_0,Z)$ vanish. Since the Poisson structures $\Pi_0$ and
$\Pi_1$ have the same transverse component, and $W$ and $Z$ are
vertical vector fields, we have $\omega_1(W,Z) = \omega_0(W,Z) =
\Pi_N^{-1}(W,Z)$. Hence $\omega_t((1-t) X_0 + tX_1,Z)=0$ as desired.

\vspace{3mm}
{\it Construction and smoothness of $\mathbb F_t$:}
\vspace{3mm}

If $X$ is a vector field on $S$ then we will denote by $X_t = (1-t) X_0 + tX_1$ the horizontal lift of $X$ to $N \times S$
via $\Gamma_t$ as above. For any two smooth vector fields $X, Y$ on $S$ and a point $q \in N \times S$, put
\begin{equation}
\bbF_t(X,Y)(q)=\omega_t(X_t,Y_t)(q) .
\end{equation}

The main point here is to check the smoothness of the function $\bbF_t(X,Y)$ defined by the above formula, in a neighborhood of the origin in $N \times S$. Denote $Z^{X}= X_0 - X_1$ and
$Z^{Y}= Y_0 - Y_1$; they are vertical vector fields. Since the Ehresmann connection $\Gamma_i$ ($i=0,1$) preserves the transverse Poisson structures, the vector fields $\hat{X_i}$ and $\hat{Y_i}$ preserve the transverse Poisson structure $\Pi_N$. Therefore the vertical vector fields $Z^{X}$ and $Z^{Y}$ also preserve the transverse Poisson structure. (They may be viewed as Poisson fields on $(N,\Pi_N)$ parametrized by $S$).

 We can write $X_t =  X_0-tZ^X=
X_1+(1-t)Z^X$ and  $Y_t =  Y_0-tZ^Y=
Y_1+(1-t)Z^Y$. Recall that if $X_t$ is horizontal with respect to $\Gamma_t$ and $Z$ is vertical then $\omega_t(X_t,Z) = 0$. We have:
\begin{equation*}
\begin{array}{ccl}
&  & \bbF_t(X,Y) \\
&= & t\omega_1(X_1+(1-t)Z^X,Y_1+(1-t)Z^Y)+ (1-t)\omega_0(X_0-tZ^X,Y_0-tZ^Y) \\
&= &  t\omega_1(X_1,Y_1)+ (1-t)\omega_0(X_0,Y_0)+ \\
&  & + t(1-t)^2\omega_1(Z^X, Z^Y)+t^2(1-t)\omega_0(Z^X, Z^Y) \\
& = & t\omega_1(X_1,Y_1)+ (1-t)\omega_0(X_0,Y_0)+ t(1-t)\Pi_N^{-1}(Z^X, Z^Y)
\end{array}
\end{equation*}

By our tameness hypothesis, $\Pi_N^{-1}(Z^X, Z^Y)$ is smooth, and so $\bbF_t(X,Y)$ is smooth (and depends smothly on $t$).

Remark that $\bbF_t$ coincides with $\bbF_0$ and $\bbF_1$ at the origin $p$. Since $\bbF_0$ is nodegenerate, $\bbF_t$ is also nondegenerate in a neighborhood of $p$ in $N \times S$. Since the form $\omega_t$ used in the construction of $(\Pi_N,\Gamma_t,\bbF_t)$ is closed on each symplectic leaf, the four compatibility conditions for the triple $(\Pi_N,\Gamma_t,\bbF_t)$ are automatically satisfied. Hence the triple $(\Pi_N,\Gamma_t,\bbF_t)$ corresponds to a smooth Poisson structure $\Pi_t$ in a neighborhood of $p$ in $N \times S$. Moreover, by construction, $\Pi_0 = \Pi$, $\Pi_1 = \Pi_N + \Pi_S$, and $\Pi_t$ depends smoothly on $t$. Lemma \ref{path} is proved.
\end{proof}

\subsection{End of the proof} \hfill

According to Lemma \ref{path}, we now have a smooth path of
$G$-invariant Poisson structures $\Pi_t$, where $\Pi_0$ is our
initial Poisson structure, and $\Pi_1 = \Pi_N + \Pi_S$ is the
splitted one. (The action of $G$ is already linearized, and by the
equivariant Darboux theorem we may assume that $\Pi_S$ is already
equivariantly normalized, i.e. has Darboux form). In order to
finish the proof of the theorem, it suffices to find a local
diffeomorphism of $N \times S$ which commutes with the action of
$G$ and which moves $\Pi_0$ to $\Pi_1$.

According to Weinstein's splitting theorem (or rather its parametrized version, whose proof is the same), there is a smooth family of local diffeomorphisms $\phi_t, t \in [0,1]$ such that
$\phi_t{_*}(\Pi_0)=\Pi_t$ and $\phi_0 = Id$. Note that, a-priori, $\phi_t$ does not commute with the action of $G$. Denote by $X_t$ the time-dependent vector field whose flow generates $\phi_t$, i.e.,
\begin{equation}
X_t(\phi_t(q))= {\partial \phi_t \over \partial t} (q).
\end{equation}

By derivation of the condition
\begin{equation}
 \phi_t{_*}(\Pi_0)=\Pi_t \label{eqn:isotopy0}
 \end{equation}
we get  the following equation for $X_t$:
\begin{equation}   \label{eqn:isotopy}
L_{X_t}(\Pi_t)=-\frac{d \Pi_t}{dt}
\end{equation}

Denote by $X_t^G$ the averaging of $X$ with respect to the action of $G$, i.e.,
\begin{equation}
X_t^G = \int_{G} \rho_g{_*}(X_{t}) d\mu, \label{eqn:averaging}
\end{equation}
where $d\mu$ is the probabilistic Haar measure on $G$, and $\rho_g$ denotes the action of  $g \in G$. Then $X^G_t$
is a $G$-invariant time-dependent vector field. Since $\Pi_t$ is invariant under the action of $G$, it follows from Equation (\ref{eqn:isotopy}) that we also have
\begin{equation}
L_{X^G_t}\Pi_t=-\frac{d \Pi_t}{dt}. \label{eqn:isotopy3}
\end{equation}

Denote by $\phi^G_t$ the flow $X^G_t$. Then $\phi^G_t$ commutes with the action of $G$.
Equation (\ref{eqn:isotopy3}) implies that $\phi^G_t{_*}(\Pi_0)=\Pi_t$. In particular, $\phi^G_1$ is a $G$-equivariant local diffeomorphism such that $\phi^G_1{_*}(\Pi_0)=\Pi_1 = \Pi_N + \Pi_S$. This concludes the proof of Theorem \ref{thm:main}.

\section{Equivariant linearization of Poisson structures}

\begin{thm} \label{thm:EquivLin}
Let $(P^n,\Pi)$ be a smooth Poisson manifold, $p$ a point of $P$, $2r = \rank \Pi (p)$, and $G$ a compact Lie group which acts on $P$ in such a way that the action preserves $\Pi$ and fixes the point $p$. Assume that the linear part of transverse Poisson structure of $\Pi$ at $p$ corresponds to a semisimple compact Lie algebra $\mathfrak k$. Then there is a smooth canonical local coordinate system $(x_1,y_1,\dots, x_{2r},y_{2r}, z_1,\dots, z_{n-2r})$ near $p$, in which the Poisson structure $\Pi$ can be written as
\begin{equation}
\Pi = \sum_{i=1}^r
\frac{\partial}{\partial x_i}\wedge \frac{\partial}{\partial y_i} + {1 \over 2}\sum_{i,j,k} c^{k}_{ij} z_k \frac{\partial}{\partial z_i}\wedge \frac{\partial}{\partial z_j},
\end{equation}
where $c_{ij}^k$ are structural constants of $\mathfrak k$,
and in which the action of $G$ is linear and preserves the subspaces $\{x_1 = y_1 = \hdots x_r = y_r = 0\}$ and $\{z_1 = \hdots = z_{n-2r} = 0\}$.
\end{thm}

\begin{proof}
Invoking Theorem \ref{thm:main}, we may assume that $\Pi$ is already equivariantly splitted, i.e.
$\Pi = \Pi = \sum_{i=1}^r
\frac{\partial}{\partial x_i}\wedge \frac{\partial}{\partial y_i} + \sum_{i,j} f_{ij}(z) \frac{\partial}{\partial z_i}\wedge \frac{\partial}{\partial z_j}$. It remains to linearize the transverse Poisson structure $\Pi_N = \sum_{i,j} f_{ij}(z) \frac{\partial}{\partial z_i}\wedge \frac{\partial}{\partial z_j}$ on $N$ in an equivariant way. But this last step is provided by the following results of Conn and Ginzburg:

\begin{thm}[Conn \cite{conn}]  \label{thm:conn}
Any smooth Poisson structure, which vanishes at a point and whose linear part at that point is of semisimple compact type, is locally smoothly linearizable.
\end{thm}

\begin{thm}[Ginzburg \cite{ginzburg}]
\label{thm:ginz}
Assume that a Poisson structure $\Pi$ vanishes at a point $p$ and is smoothly
linearizable near $p$. If there is an action of a compact Lie group $G$ which fixes $p$ and preserves $\Pi$, then
$\Pi$ and this action of $G$ can be linearized simultaneously.
\end{thm}

Indeed, by Theorem \ref{thm:conn}, the transverse Poisson structure $\Pi_N$ is smoothly  linearizable because its linear part is compact semisimple. As a consequence, by Theorem \ref{thm:ginz}, $\Pi_N$ can be linearized in a $G$-equivariant way.
\end{proof}

\section{Appendix}

In this appendix we will give a proof of Proposition \ref{prop:DivisionDim3}. We will assume that our linear Poisson structure corresponds to a 3-dimensional Lie algebra $\fg$ (the case of dimension 2 is similar and simpler and can be reduced from the 3-dimensional case). Recall that any 3-dimensional Lie algebra $\fg$ over $\bbR$ belongs to one of the following types:

\begin{enumerate}
\item   Solvable: $\fg = \bbR \ltimes_A \bbR^2$ where $A = \left( \begin{matrix}a & b \\ c & d \end{matrix} \right)$ is a 2-by-2 matrix, i.e. with Lie brackets $[x,y ] =   ay+bz$,  $[x, z ]  =   cy+dz$,  $[y,z] = 0$.
\item   Simple: $\mathfrak{so}(3,\bbR)$ or $\mathfrak{sl}(2,\mathbb R)$.
\end{enumerate}

We will prove that any vector field $X$ tangent to the symplectic
foliation of $\fg^*$ (i.e. the foliation by coadjoint orbits on $\fg^*$)
can be expressed as a smooth combination of the Hamiltonian vector fields $X_x$, $X_y$
and $X_z$, where $(x,y,z)$ is a basis of $\fg$.

 Let us first consider the
case when $\fg = \bbR \ltimes_A \bbR^2$. In this case, our linear
Poisson structure $\Pi$ can be written as:
\begin{equation}
\Pi=\frac{\partial}{\partial x}\wedge ((ay+bz)\frac{\partial}{\partial
y}+(cy+dz)\frac{\partial}{\partial z}) .
\end{equation}

We distinguish two subcases.

1)  The matrix $A$ has  non-zero determinant.

  A vector field tangent to the symplectic foliation can be written as
$Z=f\frac{\partial}{\partial x}+g ((ay+bz)\frac{\partial}{\partial
y}+(cy+dz)\frac{\partial}{\partial z})$ where the function $f$ has
vanishes  for $(ay+bz,cy+dz)=(0,0)$.

 Since the mapping $(x,y,z)\mapsto
(x,ay+bz, cy+dz)$ defines new smooth coordinates,   we may write
$f=(ay+bz) f_1+(cy+dz) f_2$ for smooth functions $f_1$ and $f_2$.

Finally we obtain $Z=f_1 X_y+f_2 X_z-g X_x$ for smooth functions
$f_1,f_2$ and $g$ as desired.

2) The matrix $A$ has determinant zero.
 In the case $a=b=c=d=0$, the
Lie algebra considered is abelian and the Poisson structure is
trivial so in this case there is nothing to prove.

In the nontrivial subcase we may write,

$$\Pi=\frac{\partial}{\partial x}\wedge (B\frac{\partial}{\partial
y}+\lambda B\frac{\partial}{\partial z})$$

\noindent being $B$ a linear function in $y$ and $z$. After a linear
change we may assume that, $\Pi=\overline B \frac{\partial}{\partial
x}\wedge \frac{\partial}{\partial \overline y}$.

A vector field tangent to the symplectic foliation  is of the form
$Z=f\frac{\partial}{\partial x}+ g\frac{\partial}{\partial \overline
y}$ where the functions $f$ and $g$ vanish when $\overline B=0$.
Since $\overline B$ is a non-trivial linear function in $\overline
y$ and $z$, we may write $f=\overline B f_1$ and $g=\overline B
g_1$. Therefore we may write $Z= f_1X_x+g_1 X_{\overline y}$.

Consider now the case when $\fg$ is simple. We will use the
following lemma  which is a smooth version of de Rham's division
lemma due to Moussu \cite{moussu}:

\begin{lemma} \label{lemma:division} Let $\alpha$ be a smooth (or analytic) 1-form on a
neighbourhood of the origin in $\mathbb R^n$ for which the origin
is an algebraically isolated singularity, then for any $p$-form
$\omega$ such that $\omega\wedge\alpha=0$ we can write the
decomposition  $\omega=\beta\wedge\alpha$ for a smooth (resp.
analytic) $(p-1)$-form $\beta$.

\end{lemma}

Denote by $\Pi$ the linear Poisson structure, it can be written as
$\Pi=x\frac{\partial}{\partial y}\wedge \frac{\partial}{\partial
z}+y\frac{\partial}{\partial z}\wedge \frac{\partial}{\partial
x}+z\frac{\partial}{\partial x}\wedge \frac{\partial}{\partial y}$
(in the case of $\mathfrak{so}(3,\mathbb K)$) or as
$\Pi=z\frac{\partial}{\partial x}\wedge \frac{\partial}{\partial
y}+x\frac{\partial}{\partial z}\wedge \frac{\partial}{\partial
x}+y\frac{\partial}{\partial y}\wedge \frac{\partial}{\partial z}$
(in the case of  $\mathfrak{sl}(2,\mathbb K)$).

 Let $\Omega$ be the volume form $\Omega= dx\wedge dy\wedge dz$, then the map $\Omega^{b}:
  \mathcal{V}^p(\mathfrak{g}^*)  \longrightarrow
\Omega^{3-p}(\mathfrak{g}^*)$ from the space of multivector fields
to the space of forms defined by $\Omega^{b} (A)= i_A\Omega$ is an
isomorphism.

Let $X$ be the vector field tangent to the symplectic foliation. The
condition of tangency to the symplectic foliation implies the
relation $X\wedge \Pi=0$. Under the above linear isomorphism this
condition becomes $i_X\Omega\wedge i_{\Pi}\Omega=0$. Since
$i_{\Pi}\Omega$ has isolated singularities at the origin, we can now
apply  lemma \ref{lemma:division} to write $i_X\Omega=\beta\wedge
i_{\Pi}\Omega$ for a smooth one-form $\beta$.

Finally, we make  convenient substitutions   to obtain
 $X  =  i_X\Omega \lrcorner(\frac{\partial}{\partial
x}\wedge \frac{\partial}{\partial y}\wedge \frac{\partial}{\partial
z})  = (\beta\wedge i_{\Pi}\Omega)\lrcorner(\frac{\partial}{\partial
x}\wedge \frac{\partial}{\partial y}\wedge \frac{\partial}{\partial
z})=\beta \lrcorner \Pi$. From this we conclude the proof of
proposition \ref{prop:DivisionDim3} since from this equality if
$\beta=fdx+gdy+hdz$ then $X=f X_x+ g X_y +h X_z$ as desired.

\hfill\qed

\bibliographystyle{amsplain}


\providecommand{\bysame}{\leavevmode\hbox
to3em{\hrulefill}\thinspace}

\end{document}